\documentclass[11pt,reqno]{amsart}
\usepackage{indentfirst, amssymb,amsmath,amsthm}
\newtheorem{theo}{Theorem}[section]

\newtheorem{lem}{Lemma}[section]
\newtheorem{rem}{Remark}

\newtheorem{cor}{Corollary}[section]
\newtheorem{defi}{Definition}[section]

\newcommand{\be}{\begin{equation}}
\newcommand{\ee}{\end{equation}}
\newcommand{\beas}{\begin{eqnarray*}}
\newcommand{\eeas}{\end{eqnarray*}}
\newcommand{\bea}{\begin{eqnarray}}
\newcommand{\eea}{\end{eqnarray}}

\numberwithin{equation}{section}

\begin{document}

\setlength{\unitlength}{1mm} \baselineskip .45cm
\setcounter{page}{1}
\pagenumbering{arabic}

\title[]
{Characterizations of generalized Robertson-Walker spacetimes concerning gradient solitons}
\author[  ]
{ Krishnendu De$^1$, Mohammad Nazrul Islam Khan$^2,^{*}$ and Uday Chand De$^3$}

\footnotetext {$\bf{2020\ Mathematics\ Subject\ Classification\:}.$ 83C05, 53C50, 53Z05.
\\ {Key words:  Gradient Ricci soliton, Gradient $(m,\tau)$-quasi Einstein soliton, Generalized Robertson-Walker spacetime, Nonlinear Equations, Differential Equations, Partial Differential Equations.\\
\thanks{$^{*}$ Corresponding author}
}}
\maketitle
\begin{center}
$^1$ Department of Mathematics,\\
 Kabi Sukanta Mahavidyalaya, The University of Burdwan.\\
 Bhadreswar, Hooghly, West Bengal, India.\\
 ORCID iD: https://orcid.org/0000-0001-6520-4520\\
\email{krishnendu.de@outlook.in }\\

$^2$Department of Computer Engineering, College of Computer,\\
Qassim University, Buraydah, Saudi Arabia.\\
ORCID iD: https: //orcid.org/0000-0002-9652-0355\\
\email{m.nazrul@qu.edu.sa}\\

$^3$Department of Pure Mathematics,
 University of Calcutta,\\
 West Bengal, India.\\
ORCID iD: https: //orcid.org/0000-0002-8990-4609\\
\email {uc$_{-}$de@yahoo.com}
\end{center}
\begin{abstract}
In this article, we examine gradient type Ricci solitons and $(m,\tau)$-quasi Einstein solitons in generalized Robertson-Walker ($GRW$) spacetimes. Besides, we demonstrate that in this scenario the $GRW$ spacetime presents the Robertson-Walker ($RW$) spacetime and the perfect fluid ($PF$) spacetime presents the phantom era. Consequently, we show that if a $GRW$ spacetime permits a gradient $\tau$- Einstein solitons, then 
it also represents a $PF$ spacetime under certain condition.
\end{abstract}

\section{\textsf{Introduction}}

Suppose $\mathcal{M}^{n}$ is a Lorentzian manifold of dimension $n$ and  $g$ is a Lorentzian metric of signature $(+,+,...,+,-)$. In 1995, the notion of $GRW$ spacetimes was proposed by Alias et al.\cite{alias1}. A $GRW$ spacetime is a Lorentzian manifold $\mathrm{M}^{n}$ $\left(n\geq4\right)$ which can be presented as $\mathcal{M}=-I \times  f^2  M^{*}$, in which $I$ $\subseteq$ $\mathbb{R}$ (Real numbers set), $M^{*}$ indicates the Riemannian manifold of dimension (n-1) and the smooth function $f>0$ is termed as warping function or scale factor. If $M^{*}$ is of dimension three and is of constant sectional curvature, then the above stated spacetime represents a $RW$ spacetime. A comprehensive investigation of $GRW$ spacetimes are presented in (\cite{alias2}- \cite{sanches1}).\par

\begin{defi}
For a scalar function $\psi$ and a 1-form $\omega_{k}$ (non vanishing), let the condition $\nabla_{k}u_{h} = \omega_{k}u_{h} + \psi g_{kh}$ be obeyed, the vector field $u$ is then referred to as torse-forming.
\end{defi}
The foregoing equation can be expressed as $\nabla_{X}u = \omega (X) u + \psi X$, $\omega$ being a 1-form.
The following theorem has been demonstrated by Mantica and Molinari \cite{survey}:

\begin{theo}(\cite{survey})\label{t2}
The Lorentzian manifold $\mathcal{M}^{n}$ ($n \ge 3$) is a $GRW$ spacetime iff the spacetime permits a unit torse-forming time-like vector field : $\nabla_{j}u_{k}=\psi (g_{ik}+u_{k}u_{i})$, it is also an eigenvector of the Ricci tensor.
\end{theo}

The $\mathcal{M}$ is termed as a $PF$ spacetime if for the non-vanishing Ricci tensor $S$, the spacetime fulfills
\begin{equation}\label{1.0}
S = a_1g + b_1\eta \otimes \eta,
\end{equation}
where $a_1, \, b_1$ are scalar fields and $g(U_1,\rho)=\eta(U_1)$ for any $U_1$ and $g(\rho,\rho)=-1$ in which $\rho$ stands for a unit time-like vector field of the $PF$ spacetime and $\eta$ is a 1-form. Each and every $RW$ spacetime presents a $PF$ spacetime \cite{neil}. However, in the dimension $4$, the $GRW$ spacetime presents a $PF$ spacetime iff the spacetime is $RW$ \cite{gtt}.\par

In a $PF$ spacetime the expression of the energy-momentum tensor $T$ is described as
\begin{equation}\label{1.1}
T=(\nu+p)\eta \otimes \eta+pg,
\end{equation}
 $\nu$ denotes the energy density, $p$ indicates the isotropic pressure \cite{neil}.\par

In absence of the cosmological constant in the theory of general relativity, the Einstein's field equations which is a highly nonlinear equations, is written as
\begin{equation}\label{1.2}
S-\frac{r}{2}g=k^2 T,
\end{equation}
where $k =\sqrt{8\pi \mathrm{G}}$, $\mathrm{G}$ indicates Newton's gravitational constant and the scalar curvature is denoted by $r$.\par
Using differential equations (\ref{1.1}) and (\ref{1.2}), we reveal the equation (\ref{1.0}), where
\begin{equation}
\label{1.3}
b_1=k^2 (p+\nu), \,\, a_1=\frac{k^2 (p-\nu)}{2-n}.
\end{equation}

Additionally, for a equation of state (EOS) parameter $\omega$, $\nu$ and $p$ are interconnected by the equation $p = \omega \nu $. The EOS having the shape $p = p(\nu )$ is named isentropic. According to \cite{ch1}, if $p=0$, $p=\dfrac{\nu}{3}$, and if $p+\nu=0$, then the PF-spacetime is represented the dust matter, the radiation and the dark energy era, respectively. Furthermore, it includes the phantom era when $\omega < -1$. The physical implications are discussed in (\cite{cal}-\cite{snu}).\par

A self-reinforcing wave packet named as a soliton, also called a solitary wave, maintains its formation while travelling with a constant speed. It is created when nonlinear and dispersive effects in the medium are neutralised. Gradient is a common term in mathematics and physics to describe the direction and magnitude of a force acting on a particle. 
In other disciplines, such as chemistry and engineering, the gradient is also used to demonstrate how a substance's property changes in relation to other variables.\par

Hamilton \cite{rsh2} develops the novel idea of Ricci flow. It is referred to as a Ricci flow \cite{rsh2} if the partial differential equations $\frac{\partial}{\partial t}g_{ij}(t)=-2S_{ij}$ satisfies the metric of a Lorentzian manifold $\mathcal{M}$. The Ricci solitons ($RS$) are produced by the self-similar solutions to the Ricci flow. If a metric of $\mathcal{M}$ obeys the differential equations,
\begin{equation}
\label{1.4}
\mathfrak{L}_{W_1}g+2S+2\lambda_{1} g=0,
\end{equation}
it is referred to as a $RS$  \cite{rsh1}, in which $\lambda_{1}$ indicates a real scalar. Also, $\mathfrak{L}_{W}$ stands for the Lie derivative operator and $W_{1}$ is the potential vector field. Equation (\ref {1.4}) has the subsequent form
\begin{equation}
\label{1.5}
Hess \, f+S+\lambda_1 g=0,
\end{equation}
in which the Hessian is denoted by $Hess$ and $D$ stands for the gradient operator of $g$ if  $W_1 = Df$, for a smooth function $f$. A gradient $RS$ is a metric that fulfills the partial differential equation (\ref{1.5}). The gradient $RS$ is said to have the smooth function $f$ as its potential function.\par

$RS$s have a significant impact in both physics and mathematics. In physics, metrics that obey (\ref{1.4}) are attractive and helpful. In connection to string theory, theoretical physicists have also been investigating the $RS$ equation. Friedan, who has done study on various features of $RS$s, has made the initial contribution to these studies \cite{fri}. In \cite{blaga2}, Blaga has considered PF spacetime endowed with a torse-forming vector field to study $\eta$-RSs and $\eta$-Einstein solitons ($ES$) and deduced a poison equation from the soliton equation. Chen and Desmukh have characterized $RS$s with the help of concurrent potential fields and on Euclidean hypersurfaces, under certain restriction they classify shrinking $RS$s \cite{chde}. Also in \cite{deshmukh2}, the authors investigated compact shrinking gradient $RS$s. Karaka and Ozgur have studied $RS$s of gradient type on multiply warped product manifolds \cite{rsr5} and obtained a necessary and sufficient condition for these manifolds to be gradient $RS$s. In \cite{wa}, Wang established that an almost $RS$ of gradient type on a $(k,\mu)'$ almost Kenmotsu manifold is a rigid gradient $RS$s.\par

In \cite{new1}, the authors have obtained exact solution for the fractional differential equations and these are emerging from solitons theory. In \cite{new2}, the authors have formulated plans that are useful in solving many different kinds of nonlinear partial differential equations arising in several areas of applied sciences. In \cite{new3}, to acquire soliton solutions to the nonlocal integrable equations, the authors have developed a new formulation of solutions to Riemann-Hilbert problems with the identity jump matrix. Rezazadeh has found a new soliton solutions of the complex Ginzburg-Landau equation with Kerr law nonlinearity in \cite{new4}. Here, we may mention that zero curvature equations make the link between integrable models and geometry manifest, and the Kronecker product produces new zero curvature representations from old ones \cite{new5}.\par


If there are $\lambda_{1}$, $\tau$ and $m$ $(0<m < \infty)$, three real constants which obeys the partial differential equation
\begin{equation}
\nabla^{2}f+S-\frac{1}{m}df \otimes df  =(\lambda_{1}+\tau r) g=\beta_{1} g,\label{a8}
\end{equation}
then the semi-Riemannian metric $g$ on the Lorentzian manifold $\mathcal{M}$ is known as a gradient $(m,\tau)$-quasi Einstein soliton ($QES$), where $\otimes$ denotes tensor product. If the potential function $f$ is constant, the soliton becomes trivial, which suggests that the manifold is Einstein. Additionally, the aforementioned relation turns into a gradient $\tau$-ES  when $m=\infty$. This idea was presented in \cite {cm}, and  Venkatesha et al. examined \cite {ven} $\tau$-$ES$ on almost Kenmotsu manifolds. More recently, in this same manifold we studied gradient $(m,\tau)$-$QES$ \cite{kde}.\par

Many researchers recently examined various types of solitons in $PF$ spacetimes, including  $RS$ (\cite{blaga2},\cite{dms}), gradient RSs (\cite{dms},\cite{dez}), Yamabe and gradient Yamabe solitons (\cite{dez}, \cite{de}), gradient m-QESs \cite{dez}, gradient $\eta$-ESs \cite{dms}, gradient Schouten solitons \cite{dms}, Ricci-Yamabe solitons \cite{sing}, respectively.\par

According to the information we have, there are many findings in the literature about $PF$ spacetimes with solitons, but there are just a few results in $GRW$ spacetimes. We want to fill this gap in this article and focus on characterizing the $GRW$ spacetimes that satisfy gradient $RS$ and gradient $(m,\tau)$-$QES$.\par

In \cite{survey}, it is established that a $GRW$ spacetime with divergence free Weyl tensor is a $PF$ spacetime. The foregoing result raises the question: Is the preceding result still valid if the condition divergence free Weyl tensor is substituted by a gradient Ricci soliton, or by a gradient $(m,\tau)$-$QES$? Here, we provide evidence that the answer to this question is, in fact, `yes' in both cases under certain conditions. Precisely, we prove the subsequent main theorems.\par

\begin{theo}\label{thm3.1}
If a $GRW$ spacetime admits a gradient $RS$ with $\rho f = constant$, then it becomes a $PF$ spacetime.
\end{theo}

\begin{theo}\label{thm3.2}
  \it{If} a $GRW$ spacetime permits a gradient $(m,\tau)$-$QES$ with $\beta_{1}=(n-1)\mu=$constant and $\rho f =$ constant, then it becomes a $PF$ spacetime.
\end{theo}

\vspace{.8cm}

\section{\textsf{Preliminaries}}

Let $\mathcal{M}$ be a $GRW$ spacetime and hence using Theorem \ref{t2}, we acquire
\begin{equation}\label{b1}
\nabla_{U_{1}}\rho=\psi[U_1+\eta(U_1)\rho]
\end{equation}
and
\begin{equation}\label{b2}
S(U_1,\rho)=\xi \eta(U_1),
\end{equation}
where $\psi$ is a scalar and $\xi$ is a non-zero eigenvector.
\begin{lem}\label{l1}
In a $GRW$ spacetime, we have
\begin{equation}\label{b3}
R(U_1,V_1)\rho = \mu[\eta(V_1)U_1-\eta(U_1)V_1]
\end{equation}
and
\begin{equation}\label{b4}
S(U_1,\rho) =(n-1)\mu\eta(U_1),
\end{equation}
where we choose $\mu=(\rho \psi + \psi^2)$.
\end{lem}
Proof. Differentiating covariantly equation (\ref{b1}), we obtain
\begin{eqnarray}\label{b5}
\nabla_{V_1}\nabla_{U_1} \rho &=& (V_1\psi)[U_1+\eta(U_1)\rho]\\ \nonumber
&&+\psi[\nabla_{V_1}U_1+(\nabla_{V_1}\eta(U_1))\rho+\psi(V_1+\eta(V_1)\rho)\eta(U_1)].
\end{eqnarray}
Interchanging $U_1$ and $V_1$ yields
\begin{eqnarray}\label{b6}
\nabla_{U_1}\nabla_{V_1} \rho &=& (U_1\psi)[V_1+\eta(V_1)\rho]\\ \nonumber
&&+\psi[\nabla_{U_1}V_1+(\nabla_X\eta(V_1))\rho+\psi(U_1+\eta(U_1)\rho)\eta(V_1)].
\end{eqnarray}
Also, we have
\begin{eqnarray}\label{bb6}
\nabla_{[U_{1},V_{1}]} \rho=\psi\{[U_{1},V_{1}]+\eta([U_{1},V_{1}])\rho\}.
\end{eqnarray}
Equations (\ref{b1}), (\ref{b5}), (\ref{b6}) and (\ref{bb6}) together implies
\begin{eqnarray}\label{b7}
R(U_1,V_1)\rho &=& (U_1\psi)[V_1+\eta(V_1)\rho] - (V_1\psi)[U_1+\eta(U_1)\rho]\\ \nonumber
&&+\psi^2[\eta(V_1)U_1-\eta(U_1)V_1].
\end{eqnarray}
Contracting $V_1$ from equation (\ref{b7}), we obtain
\begin{eqnarray}\label{b8}
S(U_1,\rho) &=& (2-n)(U_1\psi) + (\rho \psi)\eta(U_1)\\ \nonumber
&&+(n-1)\psi^2\eta(U_1).
\end{eqnarray}
Combining equations (\ref{b2}) and (\ref{b8}), we infer
\begin{equation}\label{b9}
\xi \eta(U_1)=(2-n)(U_1\psi) + (\rho \psi)\eta(U_1) +(n-1)\psi^2\eta(U_1).
\end{equation}
Setting $U_1=\rho$ in (\ref{b9}) entails that
\begin{equation}\label{b10}
\xi=(n-1)\mu,
\end{equation}
where $\mu=(\rho \psi + \psi^2)$.\par

From the last two equations, we acquire
\begin{equation}\label{b11}
U_1\psi = -(\rho \psi)\eta(U_1).
\end{equation}
Using equation (\ref{b11}) in equation (\ref{b7}), we get
\begin{equation}\nonumber\label{b12}
R(U_1,V_1)\rho = \mu[\eta(V_1)U_1-\eta(U_1)V_1].
\end{equation}
In view of equations (\ref{b2}) and (\ref{b10}), we provide
\begin{equation}\nonumber\label{b13}
S(U_1,\rho) = (n-1)\mu\eta(U_1).
\end{equation}
This ends the proof.\par

\begin{lem}\label{l2}
In a $GRW$ spacetime, we obtain
\begin{equation}\label{b14}
\mu \{U_{1}+\rho \eta(U_1)\}=0.
\end{equation}
\end{lem}
Proof: From equation (\ref{b3}), we get
\begin{equation}\nonumber\label{b15}
R(U_1,V_1)\rho = \mu[\eta(V_1)U_1-\eta(U_1)V_1].
\end{equation}
Now,
\begin{eqnarray}\label{b16}
(\nabla_{W_1} R)(U_1,V_1)\rho &=& \nabla_{W_1} R(U_1,V_1)\rho -R(\nabla_{W_1} U_1,V_1)\rho\\ \nonumber
&& -R(U_1,\nabla_{W_1}V_1)\rho -R(U_1,V_1)\nabla_{W_1} \rho.
\end{eqnarray}
Using equations (\ref{b1}) and (\ref{b3}) in equation (\ref{b16}) entails that
\begin{eqnarray}\nonumber\label{b17}
(\nabla_{W_1} R)(U_1,V_1)\rho &=&\{W_1\mu\}[\eta(V_1)U_1-\eta(U_1)V_1]\\ \nonumber
&&+\psi\mu[g(V_1,W_1)U_1-g(U_1,W_1)V_1]-\psi R(U_1,V_1)W_1.
\end{eqnarray}
The well-known second Bianchi identity is given by
\begin{equation}\nonumber\label{b18}
(\nabla_{W_1} R)(U_1,V_1)\rho + (\nabla_{U_1}R)(V_1,W_1)\rho + (\nabla_{V_1}R)(W_1,U_1)\rho =0.
\end{equation}
From the foregoing two equations, we infer
\begin{eqnarray}\nonumber\label{b19}
&&[\{W_1\mu\}\eta(V_1) - \{V_1\mu\}\eta(W_1)]U_1\\ \nonumber
&& + [\{U_1\mu\}\eta(W_1) - \{W_1\mu\}\eta(U_1)]V_1\\ \nonumber
&&+[\{V_1\mu\}\eta(U_1) - \{U_1\mu\}\eta(V_1)]W_1\\ \nonumber
&&-\psi[R(U_1,V_1)W_1 + R(V_1,W_1)U_1 + R(W_1,U_1)V_1]=0.
\end{eqnarray}
Putting $W_1=\rho$ in the previous equation gives
\begin{eqnarray}\label{b20}
&&[\{\rho\mu\}\eta(V_1) + \{V_1\mu\}]U_1\\ \nonumber
&& - [\{U_1\mu\} + \{\rho\mu\}\eta(U_1)]V_1\\ \nonumber
&&+[\{V_1\mu\}\eta(U_1) - \{U_1\mu\}\eta(V_1)]\rho\\ \nonumber
&&-\psi[R(U_1,V_1)\rho + R(V_1,\rho)U_1 + R(\rho,U_1)V_1]=0.
\end{eqnarray}
From equation (\ref{b3}), we get
\begin{equation}\label{b21}
R(\rho,U_1)V_1 = \mu[g(U_1,V_1)\rho-\eta(V_1)U_1]
\end{equation}
and
\begin{equation}\label{b22}
R(U_1,\rho)V_1=\mu[\eta(V_1)U_1-g(U_1,V_1)\rho].
\end{equation}
Using equations (\ref{b3}), (\ref{b21}) and (\ref{b22}) in equation (\ref{b20}) entails that
\begin{eqnarray}\label{b23}
&&\{\rho\mu\}[\eta(V_1)U_1-\eta(U_1)V_1]\\ \nonumber
&&+\{V_1\mu\}[U_1+\eta(U_1)\rho]\\ \nonumber
&&-\{U_1\mu\}[V_1+\eta(V_1)\rho]=0.
\end{eqnarray}
Contracting $V_1$ from the equation (\ref{b23}), we infer
\begin{equation}\nonumber
\mu\{U_{1}+\rho \eta(U_1)\}=0.
\end{equation}
Hence the proof is completed.
\begin{lem}\label{l3}
In a $GRW$ spacetime, we have
\begin{equation} \label{b25} g((\nabla_\rho Q)U_1-(\nabla_{U_{1}} Q)\rho,\rho)=0,
\end{equation}
in which the Ricci operator $Q$ is described by $g(QU_1,V_1) =S(U_1,V_1)$.
\end{lem}

Proof: From equation (\ref{b4}), we get
\begin{equation}\label{b26}
Q\rho =(n-1)\mu\rho.
\end{equation}
Differentiating equation (\ref{b26}), we acquire
\begin{eqnarray}\label{b27}
(\nabla_{U_1} Q)\rho &=& (n-1)\{U_1\mu\}\rho \\ \nonumber
&& +(n-1)\psi \mu[U_1+\eta(U_1)\rho]\\ \nonumber
&& -\psi QU_1 -(n-1)\psi \mu\eta(U_1)\rho.
\end{eqnarray}

Using equation (\ref{b27}) and Lemma \ref{l2}, we easily acquire the desired result.

\vspace{.6cm}
\section{Proof of the prime Theorems}

{\bf Proof of the Theorem \ref{thm3.1} :}\par

Let us suppose that a $GRW$ spacetime admit a gradient $RS$. Then the equation (\ref{1.5}) may be written as

\begin{equation}\nonumber
\label{c1}
\nabla_{U_{1}}Df=-QU_{1}-\lambda_{1} U_{1}.
\end{equation}

The foregoing equation and the following relation
\begin{equation}\nonumber
\label{c2}
R(U_{1}, V_{1})Df=\nabla_{U_{1}} \nabla_{V_{1}}Df-\nabla_{V_{1}} \nabla_{U_{1}}Df-\nabla_{[U_{1}, V_{1}]}Df
\end{equation}
give
\begin{equation}\nonumber
\label{c3}
R(U_{1}, V_{1})Df=(\nabla_{V_{1}}Q)(U_{1})-(\nabla_{U_{1}}Q)(V_{1}).
\end{equation}
Taking inner product of the previous equation with $\rho$ and making use of Lemma \ref{l3}, we acquire
\begin{equation}\label{c4}
g(R(U_{1},V_{1})D f,\rho)=0.
\end{equation}
Again, from equation (\ref{b3}) we infer
\begin{equation}\label{c5}
g(R(U_{1},V_{1})\rho, D f)=\mu[\eta(V_1)(U_1 f)-\eta(U_1)(V_1 f)].\end{equation}
Combining equations (\ref{c4}) and (\ref{c5}), we get
\begin{equation}\label{c6}
-\mu[\eta(V_1)(U_1 f)-\eta(U_1)(V_1 f)]=0.
\end{equation}
Replacing $V_{1}$ by $\rho$ in equation (\ref{c6}), we obtain
\begin{equation}\label{c7}
\mu[(U_{1} f)+\eta (U_{1})(\rho f)]=0.
\end{equation}
This entails that either $\mu= 0$, or $\mu\neq 0$.
\par
Case (i):
If $\mu= 0$, then from equation (\ref{b4}) we have $S(U_1,\rho) =0$. This reflects that the eigenvector $\xi$ is zero, which contradicts the Theorem \ref{t2}.\par

Case (ii): If $\mu\neq 0$, then from equation (\ref{c7}), we reveal
\begin{equation}\nonumber\label{c8}
[(U_1f) + (\rho f)\eta(U_1)]= 0,
\end{equation}
which implies
\begin{equation}\label{c9}
Df = -(\rho f)\rho.
\end{equation}

Differentiating the equation (\ref{c9}), we acquire
\begin{equation}\label{c10}
\nabla_{U_1} Df = -\{U_1(\rho f)\}\rho - \psi(\rho f)\{U_1+\eta(U_1)\rho\}.
\end{equation}
If we take $\rho f =c_1=$ constant, then either $c_1\neq 0$, or $c_1= 0$.\par

Case (i):
If $c_1\neq 0$, then equation (\ref{c10}) implies
\begin{equation}\label{c11}
\nabla_{U_1} Df = - c_1 \psi \{U_1+\eta(U_1)\rho\}.
\end{equation}
Using equation (\ref{c11}) in equation (\ref{c1}) yields
\begin{eqnarray}\nonumber
 Q U_1 &=&  c_1 \psi \{U_1+\eta(U_1)\rho\}- \lambda_1 U_1,
\end{eqnarray}
which implies
\begin{eqnarray}\nonumber\label{c12}
 S(U_1,V_1) &=& \{c_1 \psi-\lambda_1\}g(U_1,V_1)+c_1 \psi \eta(U_1)\eta(V_1).
\end{eqnarray}
Therefore, the spacetime under consideration is a $PF$ spacetime.\par

Case (ii):
If $c_1= 0$, then equation (\ref{c9}) gives $D f=0$. Using this in equation (\ref{c1}) yields
\begin{eqnarray}\label{cc121}
 S(U_1,V_1) &=& -\lambda_1 g(U_1,V_1).
\end{eqnarray}

We know that
\begin{eqnarray}\label{c}
(div C)(U_{1},V_{1})W_{1}&=&
\frac{n-3}{n-2}[\{(\nabla_{U_{1}}S)(V_{1},W_{1})-(\nabla_{V_{1}}S)(U_{1},W_{1})\}\nonumber\\&&
-\frac{1}{2(n-1)}\{g(V_{1},W_{1})dr (U_{1})-g(U_{1},W_{1})dr (V_{1})\}],
\end{eqnarray}
in which $C$ stands for the Weyl conformal curvature tensor. Therefore using equation (\ref{cc121}), from equation (\ref{c}) we acquire $(div C)(U_{1},V_{1})W_{1}=0.$\par
Thus, the spacetime is a $GRW$ spacetime with $divC=0$ and hence, it is a $PF$ spacetime \cite{survey}.\par
Hence the proof is finished.\par
\vspace{.3cm}

Since, in $4-$dimension, a $GRW$ spacetime is a $PF$ spacetime iff the spacetime is a $RW$ spacetime \cite{gtt}. Therefore, from the above theorem, we arrive:
\begin{cor}\label{cor1}
In $4-$dimension, for $\rho f = constant$, a $GRW$ spacetime admitting a gradient $RS$ turns into a $RW$ spacetime.
\end{cor}

\begin{rem}
For $n=4$, comparing the equations (\ref{1.0}) and (\ref{c12}), we have
\begin{equation}\nonumber\label{b13}
 a_1 g(U_1,V_1)+b_1 \eta(U_1)\eta(V_1)= \{c_1 \psi-\lambda_1\}g(U_1,V_1)+c_1 \psi \eta(U_1)\eta(V_1).
\end{equation}
Making use of equation (\ref{1.3}), the foregoing equation yields

\begin{equation}\nonumber\label{bb14}
    k^2 (3p-\nu)=-\lambda_1,
\end{equation}
which implies
\begin{equation}\nonumber\label{bb15}
    (3p-\nu)=c\,\, (say).
\end{equation}
Hence, the $PF$ spacetime satisfies the EOS $\nu= -3p +$ constant.\par
If $c=0$, the above equation yields
\begin{equation}\nonumber\label{bb4}
    \omega =\frac{p}{\nu}=-\frac{1}{3},
\end{equation}
which entails that the $PF$ spacetime presents the phantom era \cite{cal1}.\par

\end{rem}

\vspace{.6cm}

{\bf Proof of the Theorem \ref{thm3.2} :}\par
Let the $GRW$ spacetime permit a $(m,\tau)$-$QES$. Then the equation (\ref{a8}) may be expressed as

\begin{equation}\label{k2}
\nabla_{U_{1}}D f+Q U_{1}=\frac{1}{m}g(U_{1},D f)D f+\beta_{1} U_{1}.
\end{equation}
Differentiating covariantly equation (\ref{k2}), we obtain
\begin{eqnarray}\label{k3}
\nabla_{V_{1}}\nabla_{U_{1}}D f&=&-\nabla_{V_{1}}QU_{1}+ \frac{1}{m}\nabla_{V_{1}}g(U_{1},D f)D f\nonumber\\&&+\frac{1}{m}g(U_{1},D f)\nabla_{V_{1}}D f+\beta_{1} \nabla_{V_{1}}U_{1}+(V_{1}\beta_{1})U_{1}.
\end{eqnarray}
Interchanging $U_{1}$ and $V_{1}$ in the above equation, we get
\begin{eqnarray}\label{k4}
\nabla_{U_{1}}\nabla_{V_{1}}D f&=&-\nabla_{U_{1}}QV_{1}+ \frac{1}{m}\nabla_{U_{1}}g(V_{1},D f)D f\nonumber\\&&+\frac{1}{m}g(V_{1},D f)\nabla_{U_{1}}D f+\beta_{1} \nabla_{U_{1}}V_{1}+(U_{1}\beta_{1})V_{1}
\end{eqnarray}
and
\begin{eqnarray}\label{k5}
\nabla_{[U_{1},V_{1}]}D f=-Q[U_{1},V_{1}]+ \frac{1}{m}g([U_{1},V_{1}],D f)D f+\beta_{1} [U_{1},V_{1}].
\end{eqnarray}
From equations (\ref{k2})-(\ref{k5}), we have
\begin{eqnarray}\label{k6}
R(U_{1},V_{1})D f&=& (\nabla_{V_{1}}Q)U_{1}-(\nabla_{U_{1}}Q)V_{1}+\frac{\beta_{1}}{m}\{ (V_{1} f)U_{1}-(U_{1} f)V_{1}\}\nonumber\\&&
+\frac{1}{m}\{ (U_{1} f)QV_{1}-(V_{1} f)QU_{1} \}+\{(U_{1}\beta_{1})V_{1}-(V_{1}\beta_{1})U_{1}\}.
\end{eqnarray}
Taking inner product of equation (\ref{k6}) with $\rho$ and using Lemma \ref{l3}, we infer
\begin{eqnarray}\label{k7}
g(R(U_{1},V_{1})D f,\rho)&=&\frac{\beta_{1}}{m}\{ (V_{1} f)\eta (U_{1})-(U_{1} f)\eta(V_{1})\}\nonumber\\&&
+\frac{1}{m}\{ (U_{1} f)\eta(QV_{1})-(V_{1} f)\eta(QU_{1}) \}\nonumber\\&&
+\{(U_{1}\beta_{1})\eta (V_{1})-(V_{1}\beta_{1})\eta(U_{1})\}.
\end{eqnarray}
Again, from equation (\ref{b3}) we acquire
\begin{equation}\label{k8}
g(R(U_{1},V_{1})\rho, D f)=\mu[\eta(V_1)(U_1 f)-\eta(U_1)(V_1 f)].\end{equation}
Comparing equations (\ref{k7}) and (\ref{k8}), we obtain
\begin{eqnarray}\nonumber \label{k9}
-\mu[\eta(V_1)(U_1 f)-\eta(U_1)(V_1 f)]&=&\frac{\beta_{1}}{m}\{ (V_{1} f)\eta (U_{1})-(U_{1} f)\eta(V_{1})\}\nonumber\\&&
+\frac{1}{m}\{ (U_{1} f)\eta(QV_{1})-(V_{1} f)\eta(QU_{1}) \}\nonumber\\&&+\{(U_{1}\beta_{1})\eta (V_{1})-(V_{1}\beta_{1})\eta(U_{1})\}\nonumber.
\end{eqnarray}
Replacing $V_{1}$ by $\rho$ in the previous equation, we reveal
\begin{eqnarray}\label{nk9}
&&\{\mu-\frac{\beta_{1}}{m}+\frac{n-1}{m}\mu\}[(U_{1} f)+\eta (U_{1})(\rho f)]\nonumber\\&&+\{(U_{1}\beta_{1})+(\rho \beta_{1})\eta(U_{1})\}=0.
\end{eqnarray}
If we take $\beta_{1}=(n-1)\mu=$constant (non zero), then from equation (\ref{nk9}), we infer
\begin{equation}\nonumber\label{4.11}
[(U_1f) + (\rho f)\eta(U_1)]= 0,
\end{equation}
which implies
\begin{equation}\label{4.13}
Df = -(\rho f)\rho.
\end{equation}

Differentiating the equation (\ref{4.13}), we acquire
\begin{equation}\label{4.14}
\nabla_{U_1} Df = -\{U_1(\rho f)\}\rho - \psi(\rho f)\{U_1+\eta(U_1)\rho\}.
\end{equation}
If we take $\rho f =c_1=$ constant, then either $c_1\neq 0$, or $c_1= 0$.\par

Case (i):
If $c_1\neq 0$, then equation (\ref{4.14}) implies
\begin{equation}\label{4.15}
\nabla_{U_1} Df = - c_1 \psi \{U_1+\eta(U_1)\rho\}.
\end{equation}
Using equation (\ref{4.15}) in equation (\ref{k2}) gives
\begin{eqnarray}\nonumber
 Q U_1 &=&  c_1 \psi \{U_1+\eta(U_1)\rho\}+ \beta_{1} U_1,
\end{eqnarray}
which implies
\begin{eqnarray}\nonumber\label{cc13}
 S(U_1,V_1) &=& \{c_1 \psi+\beta_{1}\}g(U_1,V_1)+c_1 \psi \eta(U_1)\eta(V_1).
\end{eqnarray}
Hence, the spacetime taking into account is a $PF$ spacetime.\par
Case (ii):
If $c_1 = 0$, then equation (\ref{4.13}) yields $D f=0$. Using this in equation (\ref{k2}) gives
\begin{eqnarray}\nonumber
 S(U_1,V_1) &=& \beta_{1} g(U_1,V_1).
\end{eqnarray}
Using the foregoing equation in (\ref{c}), we get $div C=0$. Therefore, it is a $GRW$ spacetime with $divC=0$ and hence, it is a $PF$ spacetime \cite{survey}.\par
This ends the proof.\par
\vspace{.3cm}

It is known that when $m=\infty$, a gradient $(m,\tau)$-$QES$ produces a gradient $\tau$-$ES$. In (\ref{nk9}), we put $m=\infty$ and easily acquire the equation (\ref{4.13}). Therefore, we have:

\begin{cor}\label{cor2}
If a $GRW$ spacetime permits a gradient $\tau$- $ES$, then the gradient of the $\tau$- $ES$ potential function is pointwise collinear with the potential vector field $\rho$.
\end{cor}

Similarly, as corollary \ref{cor1} we acquire:
\begin{cor}
In $4-$dimension, a $GRW$ spacetime admitting a gradient $(m,\tau)$-$QES$ with $\beta_{1}=(n-1)\mu=$ constant and $\rho f = constant$ turns into a $RW$ spacetime.
\end{cor}

Similarly, as above we can state:

\begin{cor}\label{cor4}
If a $GRW$ spacetime admits a gradient $\tau$- $ES$ with $\rho f =$ constant, then it becomes a $PF$ spacetime.
\end{cor}

\begin{rem}
For $n=4$, comparing the equations (\ref{1.0}) and (\ref{cc13}), we have
\begin{equation}\nonumber\label{b13}
 a_1 g(U_1,V_1)+b_1 \eta(U_1)\eta(V_1)= \{c_1 \psi+\beta_1\}g(U_1,V_1)+c_1 \psi \eta(U_1)\eta(V_1).
\end{equation}
Using (\ref{1.3}), the previous equation gives

\begin{equation}\nonumber\label{bb14}
    k^2 (3p-\nu)=\beta_1,
\end{equation}
which implies
\begin{equation}\nonumber\label{bb15}
    (3p-\nu)=c\,\, (say).
\end{equation}
Therefore, the $PF$ spacetime admitting a gradient $(m,\tau)$-$QES$ obeys the EOS $\nu= -3p +$ constant.\par
If $c=0$, the above equation yields
\begin{equation}\nonumber\label{bb4}
    \omega =\frac{p}{\nu}=-\frac{1}{3},
\end{equation}
which implies that the $PF$ spacetime represents the phantom era \cite{cal1}.
\end{rem}
\vspace{.6cm}
\section{Discussion}

The stage of the physical world's current modelling is spacetime, which is a torsion less, time oriented Lorentzian manifold. Albert Einstein first proposed the idea of general relativity theory in 1915, in which the matter content of the universe is stated by picking the suitable energy momentum tensor and is accepted to act like a perfect fluid spacetime in the cosmological models. $GRW$ spacetimes, where large scale cosmology is staged, are a natural and extensive extension of $RW$ spacetimes.\par


This article will be read not only by readers working in this field, but also by researchers from other engineering disciplines. In future, other researchers or we,  will investigate others solitons in general relativity theory and cosmology.

\section*{Acknowledgment}
We appreciate the Editor's diligent examination of the manuscript and the anonymous referees' insightful criticism which enhanced the paper's quality.

\section*{Declarations}

\subsection*{Author contribution statement}

K. De: Conceived and designed the experiments; Performed the experiments; Wrote the paper.\\
M. N. Khan:  Performed the experiments;  Wrote the paper.\\
U.C. De:  Performed the experiments; Wrote the paper.
\subsection*{Funding }
Researcher would like to thank the Deanship of Scientific Research, Qassim University, for funding publication of this project.
\subsection*{Data availability statement}
No data was used for the research described in the article.

\subsection*{Declaration of interests statement}
The authors declare no conflict of interest.
\subsection*{Additional information}
No additional information is available for this paper.


\begin{thebibliography}{99}
\bibitem{alias1} L. Alias, A. Romero and M. S\'{a}nchez, {\it Uniqueness of complete spacelike hypersurfaces of constant mean curvature in generalized Robertson-Walker spacetimes}, Gen. Relat. Gravit. {\bf  27} (1995), 71.
\bibitem{alias2}
L. Alias, A. Romero and M. S\'{a}nchez, {\it compact spacelike hypersurfaces of constant mean curvature in generalized Robertson-Walker spacetimes; in geometry and topology of submanifolds VII}, River Edge NJ, USA, World Scientific, {\bf 67}, 1995.
\bibitem{bychen} B. Y. Chen, {\it Pseudo-Riemannian Geometry, $\delta$-invariants and Applications}, World Scientific, 2011.
\bibitem{bychen1}
B. Y. Chen, {\it A simple characterization of generalized Robertson-Walker spacetimes},  Gen. Relativ. Gravit. {\bf  46}, 1833 (5 pages)  (2014).
\bibitem{survey} C. A. Mantica and L. G. Molinari,  {\it Generalized Robertson-Walker spacetimes-A survey},  Int. J. Geom. Methods Mod. Phys. {\bf 14} (2017), 1730001 (27 pages).
\bibitem{mma} C. A. Mantica and L. G. Molinari,  {\it On the Weyl and Ricci tensors of Generalized Robertson-Walker spacetimes}, J. Math. Phys. {\bf 57}, 102502 (2016), https://doi.org/10.1063/1.4965714
\bibitem{sanches1}
M. S\'{a}nchez,  {\it On the geometry of generalized Robertson-Walker spacetimes: geodesics},  Gen. Relat. Gravit. {\bf 30} (1998), 915-932.
\bibitem{neil} B. O'Neill, {\it Semi-Riemannian Geometry with Applications to Relativity}, Academic Press, New York, 1983.
\bibitem{gtt} M. Guti$\acute e$rrez and B. Olea, {\it Global decomposition of a Lorentzian manifold as a generalized Robertson-Walker space}, Differ. Geom. Appl. {\bf27} (2009), 146-156.
\bibitem{ch1}P.H. Chavanis, {\it Cosmology with a stiff matter era,} Phys. Rev. D {\bf 92}, 103004 (2015).
\bibitem{cal} R. R. Caldwell, M. Kamionkowski, N. N. Weinberg, {\it Dark Energy with $w<-1$ causes a
cosmic doomsday}, Phys. Rev. Lett. {\bf 91} (2003) 071301.
\bibitem{cal1} R. R. Caldwell, {\it A Phantom Menace? Cosmological consequences of a dark energy component
with super-negative equation of state}, Phys. Lett. B {\bf545} (2002) 23–29.
\bibitem{kdn} K. De and U.C. De, {\it Almost co-Kähler manifolds and quasi-Einstein solitons}, Chaos, Solitons and Fractals, 167, 2023, 113050.
\bibitem{snu} A. Sardar, M. N. I. Khan  and U. C. De,  {\it $\eta-\ast$-Ricci solitons and almost co-Kähler manifolds}, Mathematics,  2021, 9, 3200
\bibitem{rsh2} R. Hamilton, {\it The Ricci flow on surfaces},  Contemp. Math. {\bf 71} (1988), 237-261.
\bibitem{rsh1}R. S. Hamilton,  {\it Three-manifolds with positive Ricci curvature},  J. Differ. Geom. {\bf 17} (1982), 255-306.
\bibitem{fri} D. Friedan, {\it Non linear models in $2+\epsilon$ dimensions}, Ann. Phys. {\bf163} (1985), 318–410.
\bibitem{blaga2}A. M. Blaga, {\it Solitons and geometrical structures in a perfect fluid spacetime},  Rocky Mountain J. Math. {\bf  50} (2020), 41-53.
\bibitem{chde} B.Y. Chen and S. Deshmukh, {\it Ricci solitons and concurrent vector fields}, Balkan J. Geom. Appl., {\bf 20} (2015), 14-25.
\bibitem{deshmukh2} S. Deshmukh, H. Alodan and  H. Al-Sodais, {\it A Note on Ricci Soliton}, Balkan J. Geom. Appl. {\bf16} (2011), 48-55.
\bibitem{rsr5} F. Karaca and C. Ozgur, {\it Gradient Ricci solitons on multiply warped product manifolds}, Filomat, {\bf32} (2018), 4221-4228.
\bibitem{wa} Y. Wang, {\it Gradient Ricci almost Solitons on two classes of almost Kenmotsu Manifolds}, J. Korean Math. Soc. {\bf 23} (2016), 1101-1114.
\bibitem{new1} A. Hossein, A. R. Sheikhani and H. Rezazadeh,  {\it Exact solutions for the fractional differential equations by using the first integral method}, Nonlinear Engineering, vol. 4, no. 1, 2015, pp. 15-22. https://doi.org/10.1515/nleng-2014-0018
\bibitem{new2}  H. Rezazadeh, D. Kumar, T.A. Sulaiman and H. Bulut, {\it New complex hyperbolic and trigonometric solutions for the generalized conformable fractional Gardner equation},  Modern Physics Letters B, Vol. 33, No. 17, 1950196 (2019). https://doi.org/10.1142/S0217984919501963
\bibitem{new3}W.X. Ma, {\it Nonlocal PT-symmetric integrable equations and related Riemann-Hilbert
problems}, Partial Differential Equations in Applied Mathematics 4 (2021) 100190.
\bibitem{new4}  H. Rezazadeh, {\it New solitons solutions of the complex Ginzburg-Landau equation with Kerr law nonlinearity}, Optik (2010), https://doi.org/10.1016/j.ijleo.2018.04.026
\bibitem{new5} W.X. Ma and F.K. Guo, {\it Lax Representations and Zero-Curvature Representations by the Kronecker Product}, International Journal of Theoretical Physics, Vol. 36, No. 3, 1997.
\bibitem{cm} G. Catino and L. Mazzieri,  {\it Gradient Einstein solitons},  Nonlinear Analysis {\bf 132} (2016), 66-94.
\bibitem{ven} V. Venkatesha and H.A. Kumara, {\it gradient $\rho$-Einstein soliton on almost Kenmotsu manifolds}, Annali Dell'Universita' Di Ferrara {\bf 65} (2019), 375-388.
\bibitem{kde} K. De and U.C. De, {\it A note on gradient Solitons on two classes of almost Kenmotsu Manifolds}, Int. J. Geom. Methods Mod. Phys. {\bf 19} (2022), 2250213 (12 pages).
    DOI: 10.1142/S0219887822502139
\bibitem{dms} U.C. De, C.A. Mantica and Y.J. Suh, {\it Perfect Fluid Spacetimes and Gradient Solitons}, Filomat, {\bf{36}} (2022), 829-842.
\bibitem{dez} K. De and U.C. De, A. A. Syied, N. B. Turki  and S. Alsaeed, {\it Perfect fluid spacetimes and gradient solitons}, Journal of Nonlinear Mathematical Physics, 29 (2022), 843-858.
https://doi.org/10.1007/s44198-022-00066-5
\bibitem{de} U.C. De, S.K. Chaubey and S. Shenawy, {\it Perfect fluid spacetimes and Yamabe solitons}, J Math Phys. {\bf 62}, 032501 (2021).
 https://doi.org/10.1063/5.0033967
\bibitem{sing} J.P. Singh and M. Khatri, {\it On  Ricci-Yamabe soliton and geometrical structure in a perfect fluid spacetimes}, Afrika Mathematica, {\bf{32}} (2021), 1645-1656.



	
\end{thebibliography}
\end{document}